%% file: 33.tex
\documentclass[14pt,a4paper]{scrartcl}
\usepackage[utf8]{inputenc}
\usepackage[russian]{babel}
\usepackage{amsmath}
\usepackage{amsfonts}
\usepackage{graphicx}
\usepackage{wrapfig}
\usepackage{caption}
\usepackage{import}
\usepackage{xifthen}
\usepackage{pdfpages}
\usepackage{transparent}

\begin{document}

\begin{flushleft}
\rule{170mm}{.3pt}\\
Бронза~С.~Д., Таирова~В.~Г. Профили римановых поверхностей // Теория функций, функ. Анализ и их приложения. – 1980. – Вып.~33 – с.~12"~17
\rule{170mm}{.3pt}\\
УДК 517.535.2
\end{flushleft}

\begin{center}
\large С.~Д.~Бронза, В.~Г.~Таирова\\
\large \textbf{ПРОФИЛИ  РИМАНОВЫХ  ПОВЕРХНОСТЕЙ}
\end{center}

Рассматриваем способ изображения римановых поверхностей, позволяющий описывать не только римановы поверхности класса $F_q$ (c помощью комплексов отрезков~\cite[с.~456]{Goldberg-1970}), но и римановы поверхности более широкого класса, например, замкнутые римановы поверхности положительного рода.

Рассмотрим класс римановых поверхностей, имеющих такую же структуру, как и римановы поверхности из описанного в~\cite[гл. VII, § 4]{Goldberg-1970} класса $F_q$, но без требования односвязности. Этот класс римановых поверхностей обозначим $F_q^*$. Пусть $R \in F_q^*$ и $D=\{a_i\}_{i=1}^q$~--- множество базисных точек римановой поверхности $R$. Без ограничения общности можно считать, что $D$ состоит из конечных точек. Каждую точку $a_i \in D$ снабдим $\varepsilon_i$"~окрестностью, при этом будем полагать, что замыкания этих окрестностей не пересекаются. На границе каждой $\varepsilon_i$"~окрестности выберем по точке $b_i$ и через точки множества $\Delta=\{b_i\}_{i=1}^q$ проведем замкнутую жорданову кривую $\Lambda$ так, чтобы $\varepsilon_i$"~окрестности базисных точек лежали во внешней кривой $\Lambda$. Не уменьшая общности, можно считать, что точки множества $\Delta$ располагаются на $\Lambda$, при обходе её в положительном направлении, а порядке возрастания индексов, и $b_{q+1}=b_1$. Часть кривой, расположенную между точками $b_i$ и $b_{i+1}$, будем обозначать $\Lambda_i$, внешность кривой $\Lambda-Z^*$, а внутренность~--- $Z^\circ$ (см. рис. 1).

\begin{figure}[h!]
    \centering
    \def\svgwidth{\columnwidth}
    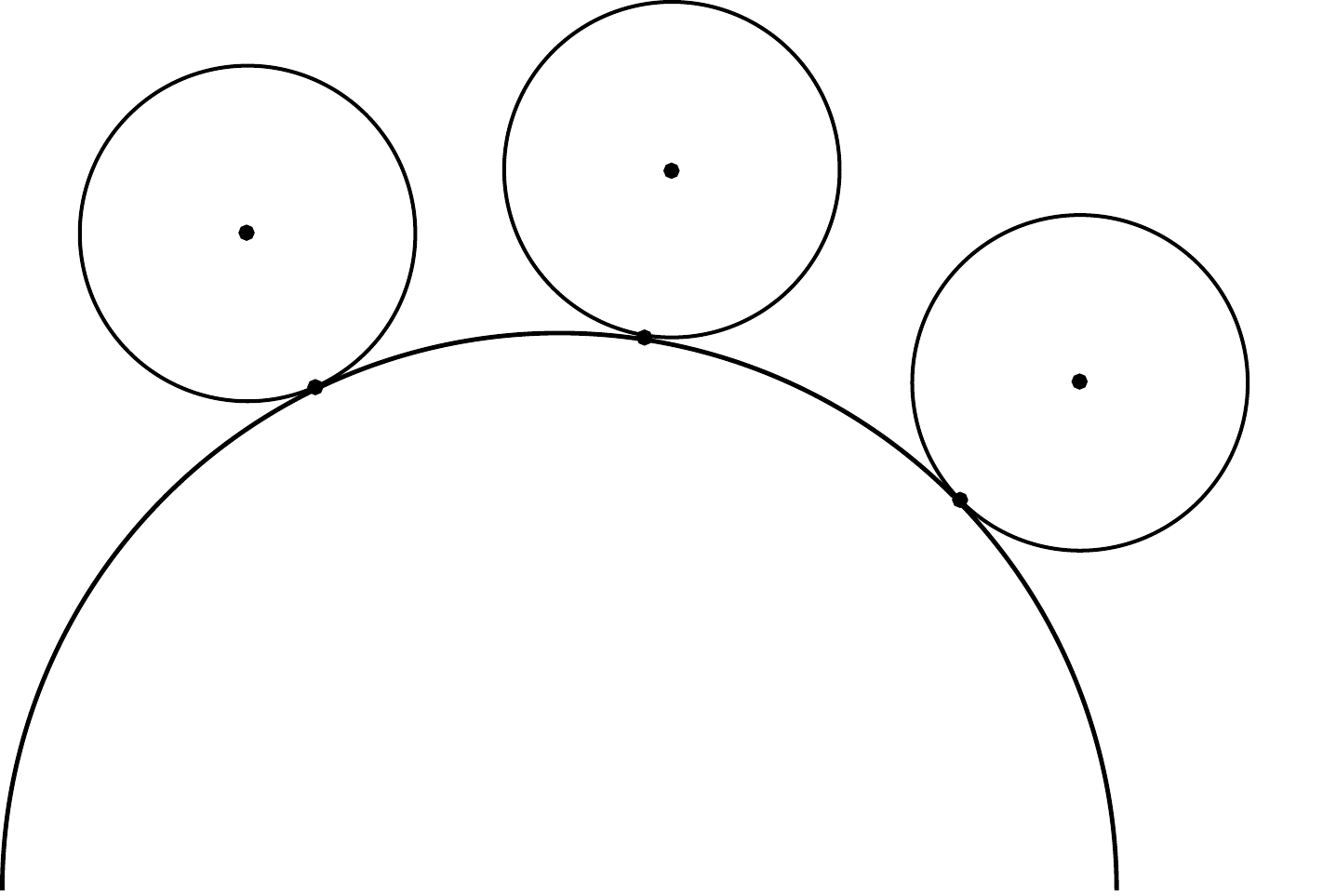

    \caption*{Рис. 1.}
    \label{fig:img_1}
\end{figure}

Рассмотрим кривую $K=\Lambda \cup (\cup_{i=1}^q \delta\varepsilon_i)$, где $\delta\varepsilon_i$~--- положительно ориентированная граница $\varepsilon_i$"~окрестности точки $a_i$. Прообраз $G$ кривой $K$ на римановой поверхности $R$, относительно операции $\pi$ проектирования ее на $\overline{C}$, можно рассмотреть как граф, множеством вершин которого является прообраз множества точек $\Delta$, а множеством ребер~--- прообраз множества $\{\Lambda_i\}_{i=1}^q$. Прообразы границ $\varepsilon_i$"~окрестностей составляют множество дуг (ориентированных ребер) графа $G$\footnote{%
Здесь и в дальнейшем используем терминологию теории графов, принятую, например в~\cite{Berger-1962}.}. Легко видеть, как граф $G$, кроме вершин, не имеет других точек самопересечения, каждой вершине графа инциденты две дуги\footnote{%
Петля графа $G$ рассматривается как две дуги~--- входящая и исходящая.} и два ребра, т.е. граф $G$~--- частично ориентированный, связный, топологический, однородный степени 4.

Элементарный цикл графа $G$, являющийся краем грани, образом которой есть $Z^\circ$, назовем $\alpha$"~циклом. Очевидно, два различных $\alpha$"~цикла не имеют общих вершин. Множество $\alpha$"~циклов графа $G$ равномощно множеству внутренних (или внешних) полулистов римановой поверхности $R$, которой соответствует граф. Перенумеруем $\alpha$"~циклы, и их упорядоченную совокупность обозначим через $\{\alpha_j\}_{j \in J}$, где $J$~--- конечное или счетное множество индексов. Точке ветвления (алгебраической или логарифмической) порядка $\lambda(\lambda \leq \infty)$, расположенной над базисной точкой $a_i$, соответствует или элементарный контур графа $G$ длины $\lambda+1$\footnote{%
Длина пути понимается в смысле теории графов, см.~\cite[с. 13]{Berger-1962}. Петлю рассматриваем как контур длины~\cite{Goldberg-1970}.} (когда $\lambda < \infty$), или бесконечный элементарный путь (когда $\lambda=\infty$), вершинами которых являются прообразы точек $b_i$ на римановой поверхности $R$. Такие контуры и пути назовем, соответственно, $\mu$"~контурами и $\mu$"~путями.

Совокупность $\mu$"~контуров и $\mu$"~путей графа $G$, соответствующих всем точкам римановой поверхности $R$, лежащим над базисной точкой $a_i$, обозначим через $\{\mu\}_{ai}$.

Придадим графу $G$ другой вид. Отметим, что каждая вершина графа $G$ одновременно принадлежит одному и только одному $\alpha$"~циклу и одной и только одной совокупности $\{\mu\}_{ai}$. Поэтому каждой вершине графа $G$ можно сопоставить упорядоченную пару чисел $(i, j)$, где $i$ указывает на принадлежность вершины множеству $\{\mu\}_{ai}$, а j~--- на ее принадлежность $a_j$"~циклу. Далее, множеству $\alpha$"~циклов $\{\alpha_j\}_{j \in J}$ поставим в соответствие множество $P$ прямых, параллельных прямой $p$; множеству $\{b_i\}_{i=1}^q$~--- множество $\{p_i\}_{i=1}^q$ точек прямой $p$, расположив его на $p$ в порядке возрастания индексов. Паре $(i, j)$ поставим в соответствие точку прямой на множества $P$, соответствующей $a_j$"~циклу, ортогонально проектирующуюся в точку $p_i$. Для этой точки сохраним обозначение $(i, j)$. Полученное множество точек рассмотрим как совокупность вершин графа $\text{П}$; смежность вершин графа $\text{П}$ определяется смежностью соответствующих им вершин графа $G$. Очевидно, что в графе $\text{П}$ вершины вида $(i, j)$ $\{i+1, j\}$ связаны ребрами, и для любой вершины $(i, j)$ найдется вершина $(i, k)$ графа $\text{П}$ такая, что обе эти вершины инцидентны одной дуге графа $\text{П}$. Если $j=k$, мы имеем петлю графа $\text{П}$. Граф $\text{П}$ назовем профилем римановой поверхности $R$ (в дальнейшем точки $p_i$ обозначим через $a_i$). На рис. 2 приведены примеры профилей римановой поверхности рода 1 с четырьмя базисными точками и римановой поверхности ${\it Arcsin}$($\omega$).

\begin{wrapfigure}{l}{0.45\textwidth}
    \centering
    \def\svgwidth{0.4\textwidth}
    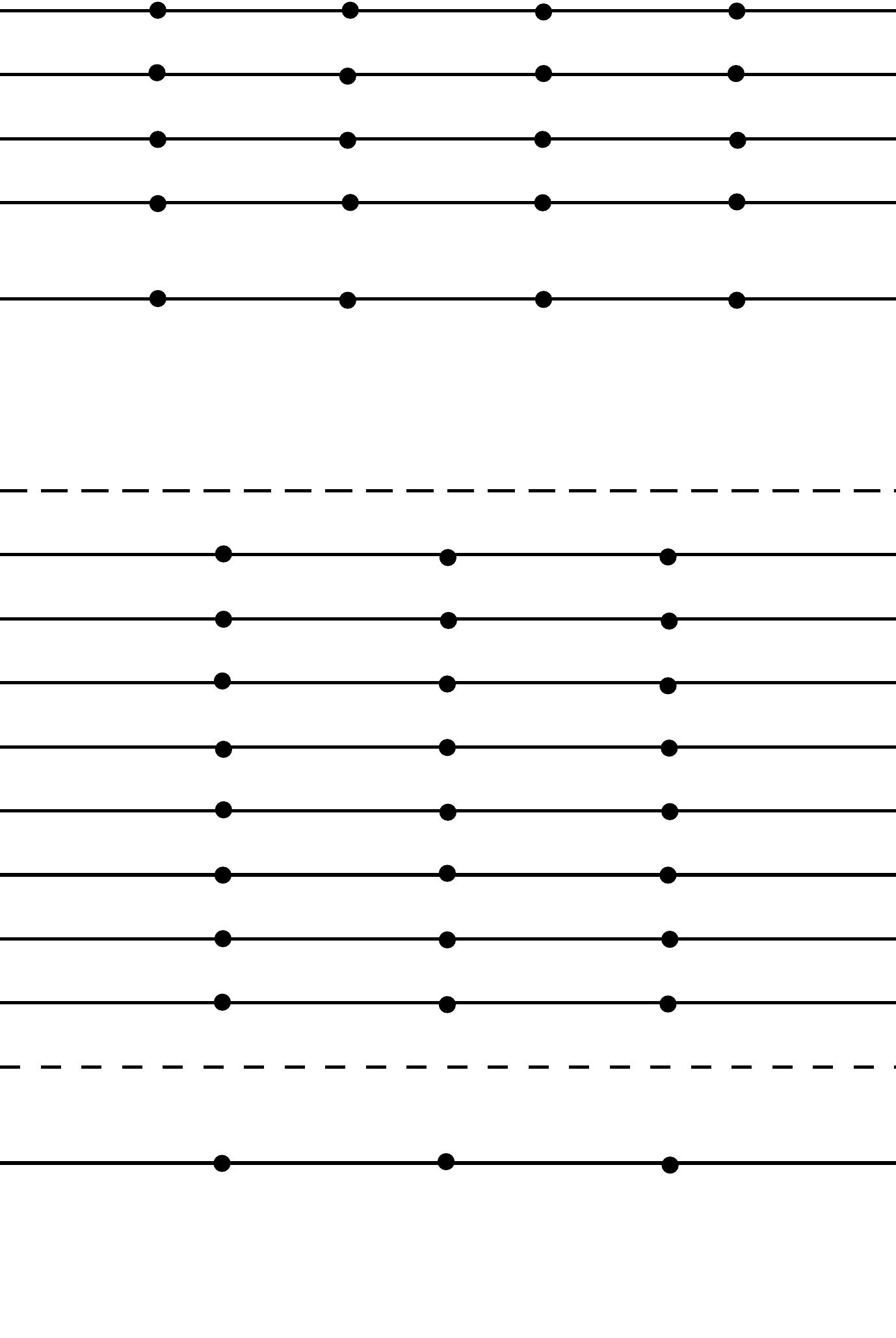

    \caption*{Рис. 2.}
    \label{fig:img_2}
\end{wrapfigure}

Пусть $\text{П}’$~--- связанный, однородный степени 4, частично ориентированный граф, каждая вершина которого инцидентна двум дугам и двум вершинам (петля рассматривается как две дуги). И пусть граф $\text{П}’$ имеет изображение, при котором его вершины расположены на конечном или на счетном множестве $P$ прямых, параллельных фиксированной прямой $p$, а множество ортогональных проекций вершин графа $\text{П}’$ на прямую $p$ конечно, при этом, если точка $p’ \in p$~--- ортогональная проекция одной из вершин графа $\text{П}’$, то все точки множества прямых $P$, имеющих проекцию $p’$, являются вершинами графа $\text{П}’$. (Прямые из $P$ отметим целочисленными индексами из некоторого множества $J$, а проекции вершин на $p$ снабдим индексами от 1 до $q$ в порядке обхода $p$). Ребрами графа $\text{П}’$ служат отрезки прямых множества $P$, на которые они разбиваются вершинами графа $\text{П}’$\footnote{%
За одно ребро засчитываем на каждой прямой из отрезок, содержащий бесконечно удаленную точку.}, а вершины, инцидентные одной дуге графа $\text{П}’$, имеет одну и туже проекцию на прямую $p$. Граф $\text{П}’$, обладающий перечисленными выше свойствами, назовем графом типа профиля. Очевидно, что профиль римановой поверхности является графом типа профиля, но не всякий граф типа профиля, есть профиль некоторой римановой поверхности.

Выясним необходимые и достаточные условия, при которых граф типа профиля является профилем некоторой римановой поверхности класса $F_q^*$. Через $(i, j)$ обозначим вершину $\text{П}’$, лежащую на прямой из $P$ с индексом $j$ и проектирующую в точку на $p$ с индексом $i$, а через~--- множество $\text{Г}$ определенных ниже путей $\text{П}’$.

О п р е д е л е н и е\quad1. Частично ориентированный цикл графа $\text{П}’$ назовем путем, если его ребра и дуги чередуются, все дуги имеют одинаковую ориентацию, при удалении петель цикл становиться элементарным\footnote{%
Петлю, входящую в путь, рассматриваем как одну дугу.}.

Вершины каждого пути $\gamma \in \text{Г}$ располагаются на нем в последовательности $(1, i_0)$, $(1, i_1)$, $(2, i_1)$, $(2, i_2)$, …, $(k, i_{k-1})$, $(k, i_k)$ …, $(q-1, i_{q-1})$, $(q, i_q)$, $(i_k \in J, i_0=i_k)$, и при обходе цикла первые индексы вершин не убывают.

О п р е д е л е н и е\quad2. Будем говорить, что совокупность путей образуют точное покрытие графа $\text{П}’$, если 1) любая дуга (ребро) графа $\text{П}’$ принадлежит некоторому пути; 2) два различных пути не имеют общих дуг (ребер).

{\bf Теорема}. Для того чтобы граф $\text{П}’$ типа профиля был профилем римановой поверхности класса $F_q^*$, необходимо и достаточно, чтобы существовало точное покрытие графа $\text{П}’$.

Множество путей точного покрытия имеет ту же мощность, что и множество $P$. Поэтому в случае замкнутых $n$"~листных римановых поверхностей с $q$ базисными точками число путей точного покрытия у соответствующих им профилей равно $n$, а длинна каждого пути равна $2q$. На рис. 3 изображены пути точного покрытия римановой поверхности.

Необходимость. Пусть $\text{П}$~--- профиль римановой поверхности $R$. Граница любой компоненты прообраза множества
\[
R=\overline{C} \thickspace \setminus \left[Z^\circ \cup \left(\bigcup_{i=1}^q \varepsilon_i\right) \cup K\right]
\]
при отображении $\pi: R \rightarrow C$, очевидно, удовлетворяет определению пути. Совокупность, таким образом определенных путей, образует точное покрытие профиля $\text{П}$. Действительно, границы двух различных компонент прообразов не имеют общих ребер (дуг), а множество ребер (дуг) профиля $\text{П}$ совпадает с множеством ребер (дуг) границы прообраза $B$, т.е. любое ребро (дуга) принадлежит одному пути.

\begin{figure}[h!]
    \centering
    \def\svgwidth{\columnwidth}
    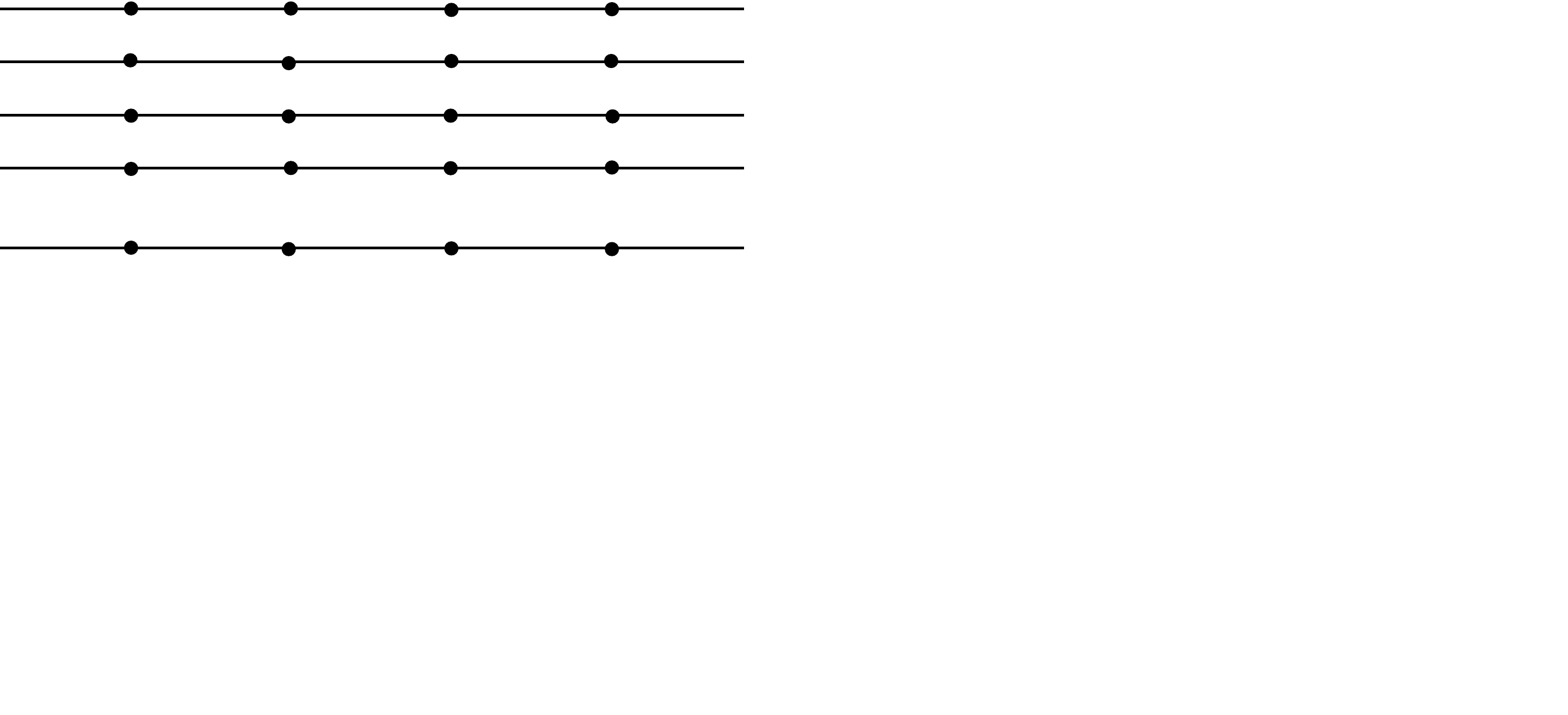

    \caption*{Рис. 3.}
    \label{fig:img_3}
\end{figure}

Достаточность. Пусть граф $\text{П}’$ имеет точное покрытие. Достаточность условий теоремы будет доказано, если графу $\text{П}’$ поставить в биективное соответствие риманову поверхность $R \in F_q^*$ и показать, что, с другой стороны, граф $\text{П}’$ является профилем этой римановой поверхности. Совокупность ребер графа $\text{П}’$, инцидентных вершинам, принадлежащим одной прямой множества P, образует цикл графа $\text{П}’$. Предположим, что множество таких циклов графа $\text{П}’$ представлено упорядоченной совокупностью $\{\beta_j\}_{j \in J}$. Возьмем на $z$"~плоскости замкнутую жорданову кривую $L$ и последовательностью $q$ точек в направлении положительного обхода разобьем ее на части $L_i$ $(i=1, … , q)$. Внутренность кривой $L$ обозначим через $Z_0$ и назовем внутренним полулистом. Внешность кривой $L$ обозначим через $Z_*$ и назовем внешним полулистом. Каждому циклу $\beta_j, j \in J$, графа $\text{П}’$ поставим в соответствие один экземпляр внутреннего полулиста $Z_0$, а каждому пути $\gamma$ из точного покрытия $\text{Г}$ графа $\text{П}’$~--- один экземпляр внешнего полулиста $Z_*$. Пусть ${\beta_i}_0$~--- некоторый цикл графа $\text{П}’$, а $l$~--- ребро этого цикла, инцидентное вершинам $(i, j_0)$, $(i+1, j_0)$. Существует один и только один путь $\gamma_0 \in \text{Г}$, имеющий с циклом ${\beta_i}_0$, общее ребро $l$. Внутренний полулист $Z_0$, соответствующий циклу ${\beta_i}_0$, и внешний полулист $Z_*$, соответствующий пути $\gamma_0$, склеим по кривой $L_i$. Такое склеивание проведем относительно ребер всех циклов $\beta_j, j \in J$, графа $\text{П}’$. Предположим, что вершина $(i, j_0)$, инцидентная ребру $l$, принадлежит контуру графа $\text{П}’$ длины $\lambda+1 (\lambda > 0)$. Вершины такого контура инцидентны $2(\lambda+1)$ ребрам, которые принадлежат ($\lambda+1$)"~му различному пути точного покрытия $\text{Г}$ и ($\lambda+1$)"~му различному циклу графа $\text{П}’$. Склеивая $\lambda+1$ внутренних полулистов, соответствующих этим путям с $\lambda+1$ внешними полулистами, соответствующими этим циклам графа $\text{П}’$, по кривым $L_i$ и $L_{i+1}$, получаем точку ветвления порядка $\lambda$ (в случае, когда $\lambda=\infty$, получаем логарифмическую точку ветвления). Если же вершина $(i, j_0)$ принадлежит петле, то инцидентные этой вершине два ребра принадлежат одному и тому же пути точного покрытия графа $\text{П}’$ и одному и тому же циклу ${\beta_i}_0$ графа $\text{П}’$.

В этом случае по кривым $L_i$ и $L_{i+1}$ склеиваются один и тот же внешний полулист с одним и тем же внутренним полулистом, в результате чего точек ветвления не возникает. Очевидно, каждое ребро внутреннего полулиста склеивается с определенным ребром некоторого внешнего полулиста, и наоборот.

Произведя указанное выше склеивание, получим некоторую риманову поверхность $R$, принадлежащую классу $F_q^*$.

Построим теперь профиль полученной римановой поверхности. В качестве базисной кривой возьмем кривую $L$. Кривую $\Lambda$ построим так, как показано на рис. 1. Пусть $Z^\circ$~--- внутренность кривой $\Lambda$, а $Z_0$~--- внутренность кривой $L$. Не уменьшая общности, можно считать, что $Z^\circ \subset Z_0$. Тогда грани графа $G$, образы которых при проектировании $\pi$ являются $Z^\circ$, принадлежат внутренним полулистам римановой поверхности $R$. Поскольку краями этих граней являются $\alpha$"~циклы графа $G$ и на каждом внутреннем полулисте располагаются по одному $\alpha$"~циклу, а сами внутренние полулисты однозначно определяются $\beta$"~циклами графа $\text{П}’$, то перенумеровав $\alpha$"~циклы графа G в соответствии с нумерацией $\beta$"~циклов, получим профиль $\text{П}$ римановой поверхности, совпадающей с графом $\text{П}’$.

Теорема доказана.

\renewcommand{\refname}{\large Список литературы:}

\begin{flushright}
Поступила 25 июня 1977 г.
\end{flushright}

\end{document}

%% file: figures/image_1.pdf_tex
\begingroup%
  \makeatletter%
  \providecommand\color[2][]{%
    \errmessage{(Inkscape) Color is used for the text in Inkscape, but the package 'color.sty' is not loaded}%
    \renewcommand\color[2][]{}%
  }%
  \providecommand\transparent[1]{%
    \errmessage{(Inkscape) Transparency is used (non-zero) for the text in Inkscape, but the package 'transparent.sty' is not loaded}%
    \renewcommand\transparent[1]{}%
  }%
  \providecommand\rotatebox[2]{#2}%
  \newcommand*\fsize{\dimexpr\f@size pt\relax}%
  \newcommand*\lineheight[1]{\fontsize{\fsize}{#1\fsize}\selectfont}%
  \ifx\svgwidth\undefined%
    \setlength{\unitlength}{409.39569505bp}%
    \ifx\svgscale\undefined%
      \relax%
    \else%
      \setlength{\unitlength}{\unitlength * \real{\svgscale}}%
    \fi%
  \else%
    \setlength{\unitlength}{\svgwidth}%
  \fi%
  \global\let\svgwidth\undefined%
  \global\let\svgscale\undefined%
  \makeatother%
  \begin{picture}(1,0.67724664)%
    \lineheight{1}%
    \setlength\tabcolsep{0pt}%
    \put(0,0){\includegraphics[width=\unitlength,page=1]{image_1.pdf}}%
    \put(0.03678507,0.03093464){\color[rgb]{0,0,0}\makebox(0,0)[lt]{\lineheight{1.25}\smash{\begin{tabular}[t]{l}$\Lambda$\end{tabular}}}}%
    \put(0.43591523,0.00502599){\color[rgb]{0,0,0}\makebox(0,0)[lt]{\lineheight{1.25}\smash{\begin{tabular}[t]{l}$Z^\circ$\end{tabular}}}}%
    \put(0.82956154,0.61598152){\color[rgb]{0,0,0}\makebox(0,0)[lt]{\lineheight{1.25}\smash{\begin{tabular}[t]{l}$Z^*$\end{tabular}}}}%
    \put(0.19063705,0.51148509){\makebox(0,0)[lt]{\lineheight{1.25}\smash{\begin{tabular}[t]{l}$a_{k+1}$\end{tabular}}}}%
    \put(0.51380873,0.55741814){\makebox(0,0)[lt]{\lineheight{1.25}\smash{\begin{tabular}[t]{l}$a_k$\end{tabular}}}}%
    \put(0.82385668,0.39665254){\makebox(0,0)[lt]{\lineheight{1.25}\smash{\begin{tabular}[t]{l}$a_{k-1}$\end{tabular}}}}%
    \put(0.23985102,0.33033936){\makebox(0,0)[lt]{\lineheight{1.25}\smash{\begin{tabular}[t]{l}$b_{k+1}$\end{tabular}}}}%
    \put(0.45885454,0.36213139){\makebox(0,0)[lt]{\lineheight{1.25}\smash{\begin{tabular}[t]{l}$b_k$\end{tabular}}}}%
    \put(0.66285841,0.26146996){\makebox(0,0)[lt]{\lineheight{1.25}\smash{\begin{tabular}[t]{l}$b_{k-1}$\end{tabular}}}}%
  \end{picture}%
\endgroup%

%% file: figures/image_2.pdf_tex
\begingroup%
  \makeatletter%
  \providecommand\color[2][]{%
    \errmessage{(Inkscape) Color is used for the text in Inkscape, but the package 'color.sty' is not loaded}%
    \renewcommand\color[2][]{}%
  }%
  \providecommand\transparent[1]{%
    \errmessage{(Inkscape) Transparency is used (non-zero) for the text in Inkscape, but the package 'transparent.sty' is not loaded}%
    \renewcommand\transparent[1]{}%
  }%
  \providecommand\rotatebox[2]{#2}%
  \newcommand*\fsize{\dimexpr\f@size pt\relax}%
  \newcommand*\lineheight[1]{\fontsize{\fsize}{#1\fsize}\selectfont}%
  \ifx\svgwidth\undefined%
    \setlength{\unitlength}{396.85042205bp}%
    \ifx\svgscale\undefined%
      \relax%
    \else%
      \setlength{\unitlength}{\unitlength * \real{\svgscale}}%
    \fi%
  \else%
    \setlength{\unitlength}{\svgwidth}%
  \fi%
  \global\let\svgwidth\undefined%
  \global\let\svgscale\undefined%
  \makeatother%
  \begin{picture}(1,1.4924419)%
    \lineheight{1}%
    \setlength\tabcolsep{0pt}%
    \put(0.17871526,0.10141639){\makebox(0,0)[lt]{\lineheight{1.25}\smash{\begin{tabular}[t]{l}$-1$\end{tabular}}}}%
    \put(0.47984915,0.1032869){\makebox(0,0)[lt]{\lineheight{1.25}\smash{\begin{tabular}[t]{l}$1$\end{tabular}}}}%
    \put(0.71855921,0.10706683){\makebox(0,0)[lt]{\lineheight{1.25}\smash{\begin{tabular}[t]{l}$\infty$\end{tabular}}}}%
    \put(0.147326,1.07548796){\makebox(0,0)[lt]{\lineheight{1.25}\smash{\begin{tabular}[t]{l}$a_1$\end{tabular}}}}%
    \put(0.36511789,1.07770084){\makebox(0,0)[lt]{\lineheight{1.25}\smash{\begin{tabular}[t]{l}$a_2$\end{tabular}}}}%
    \put(0.57702447,1.08151902){\makebox(0,0)[lt]{\lineheight{1.25}\smash{\begin{tabular}[t]{l}$a_3$\end{tabular}}}}%
    \put(0.79465806,1.08151899){\makebox(0,0)[lt]{\lineheight{1.25}\smash{\begin{tabular}[t]{l}$a_4$\end{tabular}}}}%
    \put(0,0){\includegraphics[width=\unitlength,page=1]{image_2.pdf}}%
    \put(0.47852093,0.00518488){\makebox(0,0)[lt]{\lineheight{1.25}\smash{\begin{tabular}[t]{l}$b$\end{tabular}}}}%
    \put(0.48037912,1.01181437){\makebox(0,0)[lt]{\lineheight{1.25}\smash{\begin{tabular}[t]{l}$a$\end{tabular}}}}%
    \put(0,0){\includegraphics[width=\unitlength,page=2]{image_2.pdf}}%
  \end{picture}%
\endgroup%

%% file: figures/image_3.pdf_tex
\begingroup%
  \makeatletter%
  \providecommand\color[2][]{%
    \errmessage{(Inkscape) Color is used for the text in Inkscape, but the package 'color.sty' is not loaded}%
    \renewcommand\color[2][]{}%
  }%
  \providecommand\transparent[1]{%
    \errmessage{(Inkscape) Transparency is used (non-zero) for the text in Inkscape, but the package 'transparent.sty' is not loaded}%
    \renewcommand\transparent[1]{}%
  }%
  \providecommand\rotatebox[2]{#2}%
  \newcommand*\fsize{\dimexpr\f@size pt\relax}%
  \newcommand*\lineheight[1]{\fontsize{\fsize}{#1\fsize}\selectfont}%
  \ifx\svgwidth\undefined%
    \setlength{\unitlength}{836.22050079bp}%
    \ifx\svgscale\undefined%
      \relax%
    \else%
      \setlength{\unitlength}{\unitlength * \real{\svgscale}}%
    \fi%
  \else%
    \setlength{\unitlength}{\svgwidth}%
  \fi%
  \global\let\svgwidth\undefined%
  \global\let\svgscale\undefined%
  \makeatother%
  \begin{picture}(1,0.45606218)%
    \lineheight{1}%
    \setlength\tabcolsep{0pt}%
    \put(0.06991742,0.25818408){\makebox(0,0)[lt]{\lineheight{1.25}\smash{\begin{tabular}[t]{l}$a_1$\end{tabular}}}}%
    \put(0.17327629,0.25923426){\makebox(0,0)[lt]{\lineheight{1.25}\smash{\begin{tabular}[t]{l}$a_2$\end{tabular}}}}%
    \put(0.27384213,0.26104628){\makebox(0,0)[lt]{\lineheight{1.25}\smash{\begin{tabular}[t]{l}$a_3$\end{tabular}}}}%
    \put(0.37712587,0.26104626){\makebox(0,0)[lt]{\lineheight{1.25}\smash{\begin{tabular}[t]{l}$a_4$\end{tabular}}}}%
    \put(0,0){\includegraphics[width=\unitlength,page=1]{image_3.pdf}}%
    \put(0.06991742,0.00394844){\makebox(0,0)[lt]{\lineheight{1.25}\smash{\begin{tabular}[t]{l}$a_1$\end{tabular}}}}%
    \put(0.17327627,0.00499866){\makebox(0,0)[lt]{\lineheight{1.25}\smash{\begin{tabular}[t]{l}$a_2$\end{tabular}}}}%
    \put(0.27384213,0.00681068){\makebox(0,0)[lt]{\lineheight{1.25}\smash{\begin{tabular}[t]{l}$a_3$\end{tabular}}}}%
    \put(0.37712587,0.00681063){\makebox(0,0)[lt]{\lineheight{1.25}\smash{\begin{tabular}[t]{l}$a_4$\end{tabular}}}}%
    \put(0,0){\includegraphics[width=\unitlength,page=2]{image_3.pdf}}%
    \put(0.59534114,0.25818573){\makebox(0,0)[lt]{\lineheight{1.25}\smash{\begin{tabular}[t]{l}$a_1$\end{tabular}}}}%
    \put(0.69870003,0.2592359){\makebox(0,0)[lt]{\lineheight{1.25}\smash{\begin{tabular}[t]{l}$a_2$\end{tabular}}}}%
    \put(0.79926584,0.26104792){\makebox(0,0)[lt]{\lineheight{1.25}\smash{\begin{tabular}[t]{l}$a_3$\end{tabular}}}}%
    \put(0.90254958,0.26104792){\makebox(0,0)[lt]{\lineheight{1.25}\smash{\begin{tabular}[t]{l}$a_4$\end{tabular}}}}%
    \put(0,0){\includegraphics[width=\unitlength,page=3]{image_3.pdf}}%
    \put(0.59534114,0.0039502){\makebox(0,0)[lt]{\lineheight{1.25}\smash{\begin{tabular}[t]{l}$a_1$\end{tabular}}}}%
    \put(0.69869998,0.00500031){\makebox(0,0)[lt]{\lineheight{1.25}\smash{\begin{tabular}[t]{l}$a_2$\end{tabular}}}}%
    \put(0.79926584,0.00681244){\makebox(0,0)[lt]{\lineheight{1.25}\smash{\begin{tabular}[t]{l}$a_3$\end{tabular}}}}%
    \put(0.90254958,0.00681233){\makebox(0,0)[lt]{\lineheight{1.25}\smash{\begin{tabular}[t]{l}$a_4$\end{tabular}}}}%
    \put(0,0){\includegraphics[width=\unitlength,page=4]{image_3.pdf}}%
  \end{picture}%
\endgroup%